\documentclass[11pt]{amsart}

\newtheorem{theorem}{Theorem}[section]

\newtheorem{lemma}[theorem]{Lemma}

\newtheorem{proposition}[theorem]{Proposition}
\newtheorem{definition}[theorem]{Definition}
\newtheorem{remark}[theorem]{Remark}
\newtheorem{example}[theorem]{Example}

\usepackage{amssymb}
\usepackage{amsmath, amscd}

\begin{document}
\title{Local equivalence problem for sub-Riemannian structures.}
\author{Vladimir Krouglov.}
\begin{abstract}
We solve the local equivalence problem for sub-Riemannian structures on $(2n+1)$-dimensional manifolds. We show that two sub-Riemannian structures are locally equivalent iff their corresponding canonical linear connections are equivalent. When $n=1$, these connections coincide with the generalized Tanaka-Webster connection of the corresponding contact metric structure. We show that in dimension $>5$, there may not be any contact metric manifolds associated with a given sub-Riemannian structure.
\end{abstract}
\maketitle
\section{Introduction.}
A sub-Riemannian structure is a contact structure $D$ on $M^{2n+1}$ together with a fiberwise scalar product $g$ on $D$. Sub-Riemannian structures naturally occur in different brunches of mathematics: in the study of constrained systems in classical mechanics, in optimal control, geometric measure theory and differential geometry (see \cite{M} and the references therein).

From the Darboux theorem we know that contact manifolds alone do not have local invariants (i.e. every two contact manifolds of the same dimension are locally equivalent). Sub-Riemannian structures already have local invariants. As has been shown by K.~Hughen in \cite{H}, when $n=1$ there are functions $K$ (the Webster curvature) and $\lambda$ (the eigenvalue of the torsion matrix) which do not change under the local automorphisms of a sub-Riemannian structure. Essentially the same invariants for sub-Riemannian structures on $3$-manifolds were defined by A.~Agrachev in \cite{A}. In the papers \cite{W} and \cite{NT}, similar invariants were considered in the context of $CR$-geometry and further generalized to the case of contact metric manifolds by S.~Tanno in \cite{T}. Note, that all these invariants coincide when the dimension of a manifold is three.

The main result of the present paper is a generalization of the results of K.~Hughen in \cite{H} to higher dimensions. We prove the following
\begin{theorem}
Two sub-Riemannian structures are locally equivalent if and only if their corresponding \emph{canonical connections} are locally equivalent.
\end{theorem}
The paper is organized as follows. In Section $2$ we recall basic notions and constructions from the Cartan's method of equivalence. In Section $3$ we define sub-Riemannian structures as the $G$-structures. Section $4$ is devoted to the solution of the local equivalence problem. In Section $5$ we show some relations between the sub-Riemannian geometry and contact metric geometry.

\section{G-structures on manifolds.}
In the present section we will introduce some notions and results from the theory of $G$-structures. We refer the reader to the book \cite{S} for additional details.
\subsection{G-structures.} Let $M$ be a smooth $n$-dimensional manifold and let $\mathcal{F}^\ast(M)$ denote the coframe bundle of $M$. It is a principal $GL(n, \mathbb{R})$-bundle over $M$, where $GL(n, \mathbb{R})$ action is given by the `change of basis' matrices. Let $G$ be a subgroup of $GL(n, \mathbb{R})$.
\begin{definition}
$G$-structure on $M$ is a principal $G$-subbundle of the coframe bundle of $M$.
\end{definition}
Let $\phi$ be a diffeomorphism between manifolds $M_1$ and $M_2$. The map $\phi^\ast$ induces an isomorphism of the coframe bundles $\mathcal{F}^\ast(M_1)$ and $\mathcal{F}^\ast(M_2)$.
\begin{definition}
Let $B^1$ and $B^2$ be two $G$-structures on $M_1$ and $M_2$ respectively. A diffeomorphism $\phi$ is called an isomorphism of $G$-structures if $\phi^\ast(B^2) = B^1$. We will say that two $G$-structures are isomorphic if there is an isomorphism between them.
\end{definition}
\begin{definition}
Two $G$-structures $B^1$ and $B^2$ are called locally equivalent at $(x, y)$ where $x \in M_1$ and $y \in M_2$ if there are neighborhoods $U(x)$ and $V(y)$ and a diffeomorphism $\phi: U(x) \to V(y)$ such that $\phi^\ast(B^2|_{V(y)}) = B^1|_{U(x)}$.
\end{definition}
Consider an $n$-dimensional vector space $V$ with a fixed basis $(e_1, e_2, \ldots, e_n)$. Every point $p \in \mathcal{F}^\ast(M)$ is a coframe in the tangent space $T_x(M)$ where $x = \pi(p)$.
\begin{definition}
A tautological $1$-form on $\mathcal{F}^\ast(M)$ is a $V$-valued differential $1$-form $\Theta$ on $\mathcal{F}^\ast(M)$ defined in the following way. Let $p \in \mathcal{F}^\ast(M)$ and let $X \in T_p(\mathcal{F}^\ast(M))$.
$$
\Theta(X) = \tilde{p}(d\pi(X))
$$
where $\tilde{p}: T_x M \to V$ is an evaluation map $\tilde{p}(X) = \sum_i p_i(X)e_i$.

If $B$ is a $G$-structure on $M$ then the pullback of $\Theta$ induced by inclusion $B \subset \mathcal{F}^\ast(M)$ is called a tautological form $\Theta$ on $B$.
\end{definition}
\subsection{The structure function.} For every $p$ in $B$, denote by $G_p$ the fiber of $B$ at $p$. We will call $G_p$ a vertical subspace.
Let $H_p$ be some complement to $G_p$ ( i.e. a subspace such that $T_pB = H_p \oplus G_p$). A tautological $1$-form defines an isomorphism
$$
\Theta_p : H_p \to V
$$
Consider the function $c_H \in Hom(V \wedge V, V)$ defined by
$$
c_H(u \wedge v) = d\Theta(X, Y)
$$
where $\Theta(X)=u$ and $\Theta(Y)=v$.

The Lie algebra of $G$ is a subalgebra in $\mathfrak{gl}(n, R)$ and therefore acts on $V$ by matrix multiplication. In particular, there is an embedding $\mathfrak{g} \subset Hom(V, V)$. Consider the map $\mathcal{A}: Hom(V, \mathfrak{g}) \to Hom(V \wedge V, V)$ defined by
$$
\mathcal{A}(S)(u \wedge v) = S(u)v - S(v)u
$$
The following construction is central to the theory of $G$-structures. Consider
$$
E = \frac{Hom(V \wedge V, V)}{\mathcal{A}(Hom(V, \mathfrak{g}))}
$$
and the projection $\rho: Hom(V \wedge V, V) \to E$. We call $E$ the \emph{orbit space}.
\begin{proposition}\cite{S}
For every $G$-structure $B$ there is a well-defined map called the \emph{structure function} of $B$
$$
c : B \to E
$$
such that $c(p) = \rho \circ c_{H_p}$ for every complement $H_p$ to the vertical subspace $G_p$.
\end{proposition}
\begin{theorem}\cite{S}
Let $B^1$ and $B^2$ be two $G$-structures with the structure functions $c_1$ and $c_2$ respectively. If they are locally equivalent by the diffeomorphism $\phi$, then $c_1 \circ \phi^\ast = c_2$.
\end{theorem}
An important feature of the structure function is that it is an equivariant map from $B$ to $E$. More precisely, let $\sigma$ be a linear representation of $G$ in $Hom(V \wedge V, V)$ defined by
$$
\sigma(g) S (u \wedge v) = g^{-1}S(gu \wedge gv)
$$
One may show that $\mathcal{A}(Hom(V, \mathfrak{g}))$ is a $\sigma$-invariant subspace, therefore $\sigma$ descends to a representation on $E$.
\begin{proposition}\cite{S}
For any $p \in B$ and $g \in G$ the following identity holds:
$$
c(p \cdot g) = \sigma(g)^{-1}c(p)
$$
\end{proposition}

%
\section{Sub-Riemannian structures.}
Assume that $M$ is a $(2n+1)$-dimensional manifold and $D$ is a contact structure on $M$.
\begin{definition}
A triple $\mathcal{S}=(M, D, g)$ where $g$ is a fiberwise inner product on $D$ is called a sub-Riemannian structure on $M$.
\end{definition}
Each sub-Riemannian structure defines a natural $G$-structure consisting of \emph{adapted coframes}.
\begin{definition}
We say that a coframe $(\eta_1, \eta_2, \ldots, \eta_{2n+1})$ of covectors in $T_p M$ is adapted to $\mathcal{S}$ if
\begin{enumerate}
\item{$Ker(\eta_{2n+1}) = D$.}
\item{$g = \sum_{i=1}^{2n} \eta_i \otimes \eta_i$.}
\end{enumerate}
\end{definition}
We may now consider the set $B$ of all coframes adapted to $\mathcal{S}$. Obviously, it is a principal $G$-bundle with a group
$$
G = \left \{ \begin{array}{ll} \left ( \begin{array}{ll} A & b \\ 0 & c \end{array} \right ) |& A \in O(2n, \mathbb{R}), \  b \in \mathbb{R}^{2n} \  \mbox{and} \  c \in \mathbb{R}\backslash \{0\} \end{array} \right \}
$$
Every sub-Riemannian structure $\mathcal{S}$ defines a $G$-structure $B$ of adapted coframes, and conversely -- every $G$-structure $B$ with $G$ as above gives rise to a sub-Riemannian structure $\mathcal{S}$.
\section{Local equivalence of the sub-Riemannian structures.}
We are now going to solve the local equivalence problem by applying the Cartan's method of equivalence to the $G$-structures of adapted coframes. Our argument is based on the corresponding theorem in \cite{H}.
\subsection{First Reduction.} Let $V$ be a $(2n+1)$-dimensional vector space. Fix some basis $(e_1, e_2, \ldots, e_{2n}, v)$ in $V$. We want to calculate the $\mathcal{A}$-image of $Hom(V, \mathfrak{g})$ in $Hom(V \wedge V, V)$ and to find the orbit space
$$
E = \frac{Hom(V \wedge V, V)}{\mathcal{A}(Hom(V, \mathfrak{g}))}
$$
\begin{proposition}\label{FirstReduction}
The space $E$ may be identified with $Hom(V' \wedge V', V/V)$ where $V'$ is a subspace in $V$ generated by $\{e_1, e_2, \ldots, e_{2n}\}$.
\end{proposition}
\emph{Proof:} The Lie algebra of $G$ is a set of matrices
$$
\mathfrak{g}= \left \{ \begin{array}{ll} \left ( \begin{array}{ll} A & b \\ 0 & c \end{array} \right ) |& A \in \mathfrak{o}(2n, \mathbb{R}), \  b \in \mathbb{R}^{2n} \  \mbox{and} \  c \in \mathbb{R}\backslash \{0\} \end{array} \right \}
$$
Consider the following basis for $\mathfrak{g}$:
$$
I_{pq} = e_p^\ast \otimes e_q - e_q^\ast \otimes e_p, \ \ p, q \in \overline{1,2n}
$$
$$
II_k = v^\ast \otimes e_k, \  \ k \in \overline{1,2n}.
$$
$$
III = v^\ast \otimes v
$$

The space $Hom(V, \mathfrak{g})$ is generated by vectors $e_s^\ast \otimes I_{pq}$, $e_s^\ast \otimes II_{k}$ and $e_s^\ast \otimes III$, where $e_s$ is in $\{e_1, e_2, \ldots, e_{2n}, v\}$. Calculate the $\mathcal{A}$-image of $Hom(V, \mathfrak{g})$ in $Hom(V \wedge V, V)$.
$$
\mathcal{A}(e_s^\ast \otimes II_{k})(e_i, e_j) = e_s^\ast(e_i)II_{k}(e_j) - e_s^\ast(e_j)II_{k}(e_i)
= e_s^\ast(e_i) v^\ast(e_j) e_k - e_s^\ast(e_j)v(e_i)e_k
$$
Therefore, the image of $e_s^\ast \otimes II_{k}$ in $Hom(V \wedge V, V)$ is a vector $(e_s^\ast \wedge v^\ast) \otimes e_k$. Analogously, $\mathcal{A}(e_s^\ast \otimes III)=(e_s^\ast \wedge v^\ast) \otimes v$. These vectors are zero in $E$.
$$
\mathcal{A}(e_s^\ast \otimes I_{pq})(e_i, e_j) = e_s^\ast(e_i)I_{pq}(e_j) - e_s^\ast(e_j)I_{pq}(e_i) = e_s^\ast(e_i)(e_p^\ast(e_j)e_q - e_q^\ast(e_j)e_p)
$$
$$
 - e_s^\ast(e_j)(e_p^\ast(e_i)e_q - e_q^\ast(e_i)e_p)
$$
This vector corresponds to $(e_s^\ast \wedge e^\ast_p)\otimes e_q - (e_s^\ast \wedge e^\ast_q)\otimes e_p$ in $Hom(V \wedge V, V)$. Therefore, in $E$ we have that $(e_s^\ast \wedge e^\ast_p)\otimes e_q = (e_s^\ast \wedge e^\ast_q)\otimes e_p$. We want to show that for any indices $i,j,k \in \overline{1,2n}$ the image of $e_j^\ast \wedge e_k^\ast \otimes e_i$ is equal to zero in $E$. But we may write
$$
e_i \wedge e_j \otimes e_k = e_i \wedge e_k \otimes e_j = - e_k \wedge e_i \otimes e_j = - e_k \wedge e_j \otimes e_i = e_j \wedge e_k \otimes e_i =  e_j \wedge e_i \otimes e_k
$$
which finishes the proof of the proposition.

The group $GL(2n, \mathbb{R})$ acts on the space $E$ by
$$
A(\omega \otimes [v]) (x, y) = \omega(A^{-1}x, A^{-1}y)[v]
$$
where $x$ and $y$ are in $V'$ and $[v]$ is an image of $v$ in $V/V'$.
\begin{proposition}
Assume that $\mathcal{S} = (M, D, g)$ is a sub-Riemannian structure on $M$. The structure function of $\mathcal{S}$ takes values in the open $GL(2n, \mathbb{R})$-orbit in $E$. This orbit is a trivial fibration by the orbits of the $G$-action. The stabilizer of each point is conjugate to
$$
 G_1 = \left \{ \begin{array}{ll} \left ( \begin{array}{ll} A & b \\ 0 & 1 \end{array} \right ) |& A \in U(n), \  b \in \mathbb{R}^{2n} \end{array} \right \}
$$
in $G$.
\end{proposition}
\emph{Proof:} By the Proposition \ref{FirstReduction} the space $E$ is identified with the space $V/V'$-valued 2-forms in $V'$. Consider the set of all nondegenerate 2-forms in $E$. It is an open $GL(2n, \mathbb{R})$-orbit in $E$. We need to show that if $D$ is a contact structure the structure function of $B$ takes values in this set.

Take some point $p$ in $B$. We need to prove that for every nonzero $x \in V'$ there exists some $y \in V'$ such that
$$
c(p)(x, y) \ne 0
$$
Take some complement $H_p$ to the vertical subspace at $p$ and consider the vector $\tilde{X}$ in $H_p$ such that $\Theta(\tilde{X})=x$. For every $\tilde{Y}$,
$$
c_H(p)(x, y) = d\Theta(\tilde{X}, \tilde{Y})=\sum_{i=1}^{2n} d\Theta^i(\tilde{X}, \tilde{Y})e_i + d\Theta^{2n+1}(\tilde{X}, \tilde{Y})v
$$
Consider the projection $\rho: Hom(V \wedge V, V) \to E$. The structure function may be written as
$$
c(p)(x, y) = \rho \circ c_H(p)(x, y) = d\Theta^{2n+1}(\tilde{X}, \tilde{Y})v
$$
But
$$
d\Theta^{2n+1}(\tilde{X}, \tilde{Y})v = [\tilde{X}\Theta^{2n+1}(\tilde{Y}) - \tilde{Y}\Theta^{2n+1}(\tilde{X}) - \Theta^{2n+1}([\tilde{X}, \tilde{Y}])]v
$$
If $\Theta(\tilde{Y})$ is in $V'$, the first two summands are zero. Finally, as follows from contact condition, there exists some $\tilde{Y}$ such that $d\Theta^{2n+1}(\tilde{X}, \tilde{Y}) = - \Theta^{2n+1}([\tilde{X}, \tilde{Y}]) \ne 0$.

The group $G$ acts on $E$ by
$$
\left ( \begin{array}{ll} A & b \\ 0 & c \end{array} \right ) \omega \otimes [v] = cA^{-1} \omega \otimes [v]
$$
On $E$ this action coincides with the action of the \emph{conformal orthogonal group} $CO(2n, \mathbb{R})$ induced by the action of $GL(2n, \mathbb{R})$. In particular, the $G$-stabilizer of each point is a preimage of the stabilizer by the $CO(2n, \mathbb{R})$ action under the homomorphism
$$
\left ( \begin{array}{ll} A & b \\ 0 & c \end{array} \right ) \to cA.
$$
The stabilizer of a $GL(2n, \mathbb{R})$-action is isomorphic to $Sp(2n, \mathbb{R})$ and therefore the stabilizer of a $CO(2n, \mathbb{R})$-action at the point $\omega \otimes v$ is an intersection $CO(2n, \mathbb{R}) \cap Sp(\omega)$, where $Sp(\omega)$ is a group of all linear isomorphisms which preserve the symplectic form $\omega$.

Since the determinant of a symplectic matrix is equal to one, the intersection $CO(2n, \mathbb{R}) \cap Sp(\omega) = O(2n, \mathbb{R}) \cap Sp(\omega) = U(\omega)$ -- the group of all matrices that respect the standard inner product in $V'$ (with respect to a chosen basis) and the symplectic form $\omega$.

Lets now look at  the action of a group $O(2n, \mathbb{R})$ on the $GL(2n, \mathbb{R})$-orbit. It is well known (see for example \cite{SM}) that the embedding of the $O(2n, \mathbb{R})$-orbit $O(2n, \mathbb{R})/U(n) \hookrightarrow GL(2n, \mathbb{R})/Sp(2n, \mathbb{R})$ is a homotopy equivalence. Consider the long exact sequence of homotopy groups:
$$
\ldots \to \pi_n(O(2n, \mathbb{R})/U(n)) \to \pi_n(GL(2n, \mathbb{R})/Sp(2n, \mathbb{R})) \to \pi_n(\mathcal{\mathcal{B}}) \to
$$
$$
\to \pi_{n-1}(O(2n, \mathbb{R})/U(n)) \to \ldots \to \pi_0(GL(2n, \mathbb{R})/Sp(2n, \mathbb{R})) \to \pi_0(\mathcal{\mathcal{B}}) \to 0.
$$
Since $\pi_0(O(2n, \mathbb{R})/U(n)) = \pi_0(GL(2n, \mathbb{R})/Sp(2n, \mathbb{R}))$ we may see that the base $\mathcal{B}$ of fibration by the $O(2n, \mathbb{R})$ action is contractible and hence the fibration by the $O(2n, \mathbb{R})$-orbits is trivial. We want to show that fibration by the $G$-orbits is also trivial. For this consider the following commutative diagram induced by the identity map
$$
\begin{CD}
O(2n, \mathbb{R})/U(n) @>>> GL(2n, \mathbb{R})/Sp(2n, \mathbb{R}) @>>> \mathcal{B}\\
@ViVV @VidVV @VjVV\\
CO(2n, \mathbb{R})/U(n) @>>>GL(2n, \mathbb{R})/Sp(2n, \mathbb{R}) @>>> \mathcal{B'}
\end{CD}
$$
It is easy to see that the map $i$ is a homotopy equivalence, hence as follows from the `five lemma' applied to the induced commutative diagram of the long exact sequences, the map $j$ is a homotopy equivalence which shows that $\mathcal{B}'$ is also contractible and the fibration by the $CO(2n, \mathbb{R})$-orbits is trivial. Finally, the orbits of the $CO(2n, \mathbb{R})$-action coincide with the orbits of the $G$-action hence the fibration by the $G$-orbits is trivial.

Consequently, there is a $G$-equivariant diffeomorphism 
$$CO(2n, \mathbb{R})/U(n) \times \mathcal{B}' \to GL(2n, \mathbb{R})/Sp(2n, \mathbb{R}).$$ For every $x \in CO(2n, \mathbb{R})/U(n)$, the image of $\{x\} \times \mathcal{B}'$ is some section $\textbf{s}$ with the property that each of its point has exactly the same stabilizer (i.e. a \emph{slice}). The stabilizer of a slice $\textbf{s} = [id]\times \mathcal{B}'$ is exactly $G_1$. This finishes the proof of the proposition.

Consider the set $\mathcal{S}_1 = c^{-1}(\textbf{s})$. One may show that $\mathcal{S}_1$ is a $G_1$-structure (see \cite{S} for further details). We are going to call $\mathcal{S}_1$ the \emph{first reduction} of $\mathcal{S}$.

\begin{lemma} \label{Equiv} Two sub-Riemannian structures are locally equivalent if and only if their corresponding first reductions are locally equivalent.
\end{lemma}
\emph{Proof:} If the reduced structures are locally equivalent, then extending the structure groups to $G$ would give the equivalence of the original sub-Riemannian structures. Conversely, assume that $\mathcal{S}$ is locally equivalent to $\mathcal{S}'$. Denote by $\phi$ the equivalence an let $c_1$ and $c_2$ be the structure functions of $\mathcal{S}_1$ and $\mathcal{S}_2$. Since
$$
c_1 \circ \phi^\ast = c_2
$$
we have that $c_1 \circ \phi^\ast(c_2^{-1}(\textbf{s})) = c_2(c_2^{-1}(\textbf{s})) = \textbf{s}$.
\subsection{Second reduction.}
\begin{proposition}
The orbit space $E_1$ of the $G_1$-structure $\mathcal{S}_1$ may be identified with the space $A \oplus Hom(V' \wedge V', V/V') \oplus Hom(V' \wedge V/V', V/V')$, where $A$ is some subspace in $Hom(V' \wedge V', V')$.
\end{proposition}
\emph{Proof:} The Lie algebra $\mathfrak{g}_1$ has the following basis:
$$
II_k = v^\ast \otimes e_k, \ k=\overline{1,2n}
$$
$$
A_{pq} = I_{pq} - J_0I_{pq} - I_{qp} + J_0I_{qp}, \  p, q=\overline{1, 2n}
$$
where $J_0$ is a standard complex structure of $\mathbb{R}^{2n}$.

Using the same arguments as in Proposition \ref{FirstReduction}, we may prove that
the $\mathcal{A}$-images of vectors $e_s^\ast \otimes II_k$ would be zero in $E_1$. These vectors would span the space $Hom(V' \wedge V/V', V')$. The space $Hom(V \wedge V, V)$ may be decomposed as a direct sum
$$
Hom(V \wedge V, V) = Hom(V' \wedge V', V/V') \oplus Hom(V' \wedge V/V', V')
$$
$$
 \oplus Hom(V' \wedge V', V') \oplus Hom(V' \wedge V/V', V/V')
$$
and the claim follows.
\begin{remark}
Unless $n=1$ the subspace in $Hom(V \wedge V, V)$ generated by the vectors $\mathcal{A}(A_{pq})$ would not coincide with $Hom(V' \wedge V', V')$ and $A$ would not be zero-dimensional.
\end{remark}
In order to make the second reduction we are going to use only the part of the structure function.

Denote by $\chi$ the projection
$$
E_1 \to E_2=Hom(V' \wedge V', V/V') \oplus Hom(V' \wedge V/V', V/V')
$$
and consider the function $\chi \circ c$. For every point $p$ in $\mathcal{S}_1$ its image has the form
$$
(\chi \circ c )(p) = \omega(p) \otimes v + \eta(p) \wedge v^\ast \otimes v
$$
where $\omega(p) \otimes [v] \in \textbf{s} \subset E$ is an $U(n)$-invariant $V/V'$-valued $2$-form.

The image of this vector under the $G_1$-action is
$$
\left(\begin{array}{ll}A & b \\ 0 & 1\end{array}\right)(\omega \otimes v + \eta \wedge v^\ast \otimes v) = \omega \otimes v + A^{-1}\eta \wedge v^\ast \otimes v + i_b \omega \wedge v^\ast \otimes v
$$
where $i_b \omega$ is a interior product of $b$ and $\omega$ (i.e. $i_b \omega (u) = \omega(b, u)$). Since $\omega$ is nondegenerate, for every $V/V'$-valued $1$-from $\zeta$, there is a unique solution $b$ of the system of linear equations $i_b \omega + A^{-1}\eta=\zeta$. In particular, all vectors of the form
$$
\omega \otimes v + \zeta \wedge v^\ast \otimes v
$$
lie on the same $G_1$-orbit. Therefore, every orbit of the $G_1$-action is an affine subspace $(\omega \otimes v) \oplus  Hom(V' \wedge V/V', V/V') $ and the stabilizer of each point is conjugate to $\mathbb{R}^{2n}$.

Take a section of $E_2$ of the form $\textbf{s}_1=\{\omega \otimes v + 0 \cdot \eta \wedge v^\ast \otimes v\}$. The stabilizer of each point in $\textbf{s}_1$ would be
$$
G_2 = \left \{ \begin{array}{ll} \left ( \begin{array}{ll} A & 0 \\ 0 & 1 \end{array} \right ), \ \mbox{where} \  A \in U(n) \end{array} \right \}
$$
Consider the set $\mathcal{S}_2 = (\chi \circ c)^{-1}(\textbf{s}_1)$. It will be a $G_2$-structure and using the same argument as in Lemma \ref{Equiv}, two sub-Riemannian structures are locally equivalent if and only if their \emph{second reductions} are locally equivalent.
\subsection{Prolongation.} In previous sections we studied the equivalence problem for sub-Riemannian structures via the first order structure function $c$. Following the approach of Hughen we will now use the prolongation procedure to the reduced $G_2$-structure.

Let $B_G$ be a $G$-structure. Fix a complement $C$ to $\mathcal{A}(Hom(V, \mathfrak{g}))$ in $Hom(V \wedge V, V)$. Then for each point of $B_G$ we have a distinguished class of horizontal subspaces $\mathcal{H}_p = \{H: c_H \in C\}$. Now, any horisontal subspace at a point $p \in B_{G_2}$ induces a $\mathfrak{g}$-valued $1$-form $\omega$ such that
$$
\left\{
\begin{array}{l}
\omega_H(A) = A, \mbox{\ if \ } A \in G_p \\
\omega_H(A) = 0, \mbox{\ if \ } A \in H_p
\end{array}
\right.
$$
For every $H \in \mathcal{H}_p$ the pair $(\Theta_p, \omega_H)$ is a coframe in $T_p B_{G_2}$. In particular, by this identification the set $\mathcal{H}$ defines a subset in $\mathcal{F}^\ast(B_{G_2})$.

One may show that if $H_1$ and $H_2$ both satisfy $c_{H_i} \in C$, then $S_{H_1, H_2} \in \mathfrak{g}^{(1)} = Ker(\mathcal{A})$.
\begin{definition}
We call $G^{(1)}$ the group of all linear transformations of the form $a_T, T \in \mathfrak{g}^{(1)}$ where
$$
a_T(\Theta) = \Theta
$$
$$
a_T(\omega_H) = \omega_H + T\omega_H
$$
\end{definition}
It can be shown that $\mathcal{H}$ defines a $G^{(1)}$-structure in $\mathcal{F}^\ast(B_{G})$. It is called a first prolongation of $B_{G}$.
\begin{theorem}\cite{S}
A choice of complement $C$ picks out a $G^{(1)}$-structure on each $G$-structure $B_G$. Two $G$-structures are locally equivalent if an only if their first prolongations are locally equivalent (as $G^{(1)}$-structures).
\end{theorem}
Lets now look at the prolongation of the $G_2$-structure obtained by the second reduction.
\begin{proposition}
There exists a complement $C$ to $\mathcal{A}(Hom(V, \mathfrak{g}))$ in $Hom(V \wedge V, V)$ such that the first prolongation of a $G_2$-structure is a $\{e\}$-structure induced by some $G_2$-connection on $\mathcal{S}_2$.
\end{proposition}
\emph{Proof:} First observe that every first prolongation of the $G_2$-structure is a $\{e\}$-structure. The result follows either directly from the calculation of $Ker(\mathcal{A})$ or from the fact that $G_2 \subset O(2n+1, \mathbb{R})$ and $\mathfrak{o}(N, \mathbb{R})^{(1)}=0$ for every $N$.
It follows that the set $\mathcal{H}$ defines a distribution which is transverse to a vertical subspace. For it to be a $G_2$-connection it is sufficient to find $C$ which would be invariant under the action of $G_2$ on $Hom(V \wedge V, V)$, since
$$
c_{dR_g(H)}(p \cdot g) = \sigma^{-1}(g) c_H(p) \in C.
$$
We can easily find such $C$, taking for example the orthogonal complement to $\mathcal{A}(Hom(V, \mathfrak{g}))$ with respect to some $G_2$-invariant scalar product on $Hom(V \wedge V, V)$.
\begin{definition}
We say that a $G_2$-connection on $\mathcal{S}_2$ is \emph{canonical} if it may be obtained from some $G_2$-invariant $C$ by the above procedure.
\end{definition}
We may now state our main theorem:
\begin{theorem}
Two sub-Riemannian structures are locally equivalent if and only if their canonical connections which correspond to some $C$ are locally equivalent as linear connections.
\end{theorem}
\begin{remark}
The situation described above differs significantly with the situation in Riemannian and CR-geometry, where there exists a \emph{unique} canonical connection - the Levi-Civita and Webster connection correspondingly. Unless $n=1$, apriori there may exist a lot of different canonical connections associated with a given sub-Riemannian structure.
\end{remark}
\section{Sub-Riemannian structures and contact Riemannian geometry.}
Let $(M, D)$ be a $(2n+1)$-dimensional contact manifold and let $\eta$ be some fixed contact $1$-form for $D$. There is a unique vector field $\xi$ (called the Reeb vector field of $\eta$) such that
$$
\eta(\xi) = 1, \ i_\xi d\eta = 0.
$$
If $g$ is a Riemannian metric on $M$ and $\phi$ is a $(1,1)$-tensor field such that
$$
g(\xi, X) = \eta(X), \ g(X, \phi Y) = d\eta(X, Y), \ \phi\phi X = -X + \eta(X)\xi
$$
the metric is called the metric \emph{associated} with $\eta$ and the tuple $(M, \eta, \xi, g, \phi)$ - a \emph{contact Riemannian manifold}.

To every second reduction of the structure group $G$ we may canonically associate a Riemannian metric on $M$. For this, consider the pullback by inclusion of the form $\Theta^{2n+1}$ to $B_{G_2}$. This form is $G_2$-invariant and therefore defines some $1$-form $\eta$ on $M$. It is clear that $\eta$ is a contact form, since for every $X \in D$
$$
\eta(X) = \Theta^{2n+1}(\tilde{X}) = \tilde{p}^{2n+1}(d\pi(\tilde{X})) = \eta^{2n+1}(X) = 0
$$
where $\tilde{X}$ is a lift of $X$ at any point $p = (\eta^1, \ldots, \eta^{2n}, \eta^{2n+1})$.

Denote by $\xi$ the Reeb vector field of $\eta$ and define a Riemannian metric $\tilde{g}$ on $M$ by
$$
\left\{
\begin{array}{l}
\tilde{g}(X, Y) = g(X, Y), \mbox{for all} \ X, Y \in D \\
\tilde{g}(\xi, \xi) = 1 \\
\tilde{g}(\xi, X) = 0, \mbox{for all}\  X \  \in D
\end{array}
\right.
$$
A $2$-form $d\eta$ uniquely defines the operator $\phi$ on $D$ by
$$
\tilde{g}(X, \phi Y) = d\eta(X, Y), \mbox{for every} \ X, Y \in TM.
$$
We extend it to the whole tangent space by setting $\phi \xi = 0$.

\begin{proposition}
If $\dim(M)=3$ then for every sub-Riemannian structure $\mathcal{S}$ on $M$ there is a canonical Riemannian metric associated with $g$.
\end{proposition}
\emph{Proof:} Let $(D, g)$ be a sub-Riemannian structure on $M$. Fix a basis $(e_1, e_2, v)$ in $V$. The space $E = \frac{Hom(V \wedge V, V)}{\mathcal{A}(Hom(V, \mathfrak{g}))}$ is $1$-dimensional and is generated by the vector $e_1^\ast \wedge e_2^\ast \otimes v$ and the group $G$ acts on $E$ by
$$
\left(\begin{array}{ll} A & b \\ 0 & c \end{array}\right) k \cdot e_1^\ast \wedge e_2^\ast \otimes v = ck \det(A) \cdot e_1^\ast \wedge e_2^\ast \otimes v
$$
In particular, the action of $G$ has two orbits: $0$ and $\mathcal{O} = E \smallsetminus 0$ and when $D$ is a contact structure the structure function of $\mathcal{S}$ takes values in $\mathcal{O}$. If we consider the point $\omega = e_1^\ast \wedge e_2^\ast \otimes v \in E$, the first reduction of the structure group traces out a unique $1$-form $\eta$ such that $Ker(\eta) = D$. It is easy to see that the second reduction defines an $SO(2)$-structure of coframes $(\eta, \eta_1, \eta_2)$ that satisfy
$$ \left\{
\begin{array}{l}
d\eta = \eta_1 \wedge \eta_2 \\
g = \eta_1^2 + \eta_2^2
\end{array}
\right.
$$
Obviously, the operator $\phi$ that corresponds to $d\eta$ is an \emph{almost complex structure} on $D$ and the metric $\tilde{g}$ is associated with $\eta$.

As the following example shows, when $n \ge 5$ there might be no metrics associated with a sub-Riemannian structure.
\begin{example}[Sub-Riemannian structure with no associated metric.]
\end{example}
Consider the space $\mathbb{R}^5$ and a sub-Riemannian structure $(D, g)$
$$D = Ker(dz + x_1 dy_1 + x_2 dy_2)$$
$$
g = p ~dx_1^2 + q ~dy_1^2 + r ~dx_2^2 + s ~dy_2^2.
$$
Denote by $\alpha = dz + x_1 dy_1 + x_2 dy_2$. Every conact $1$-form $\eta$ such that $Ker(\eta) = D$ has a form $f\alpha$ for some function $f$.

Therefore, the differential
$$
d\eta = df \wedge \alpha + f (dx_1 \wedge dy_1 + dx_2 \wedge dy_2)
$$
when restricted to $D$ is simply $f (dx_1 \wedge dy_1 + dx_2 \wedge dy_2)$.

Consider the following local frame in $D$:
$$
\left\{
\begin{array}{l}
e_1 = \partial/\partial x_1 \\
e_2 = \partial/\partial x_2 \\
e_3 = x_1 \partial/\partial z - \partial/\partial y_1 \\
e_4 = x_2 \partial/\partial z - \partial/\partial y_2
\end{array}
\right.
$$

With respect to this basis the operator that is adjoint to $d\eta$ acts as follows:
$$
\left\{
\begin{array}{l}
\phi e_1 = \frac{f}{r} e_3 \\
\phi e_2 = \frac{f}{s} e_4 \\
\phi e_3 = -f/p e_1 \\
\phi e_4 = -f/q e_2
\end{array}
\right.
$$
and one may always find real numbers $p, q, r, s$ such that $\tilde{g}$ is not associated with $\eta$ for any $f$.


\begin{thebibliography}{LLL}
\bibitem{A}
\emph{A.A.~Agrachev}, Exponential mappings for contact sub-Riemannian structures, \emph{J.Dynamical Systems and Control}, 2(3): 321--358, 1996
\bibitem{H}
\emph{K.~Hughen}, The geometry of sub-Riemannian three-manifolds, \emph{PhD Thesis.}, Duke University, 1996
\bibitem{SM}
\emph{D.~McDuff, D.~Salamon}, Introduction to Symplectic Topology, Oxford University Press, 1998, 512 p.
\bibitem{M}
\emph{R.~Montgomery}, A tour of sub-Riemannian geometries -- their geodesics and applications, \emph{Math. Surveys and Monographs}, AMS, Providence, RI, 2002
\bibitem{S}
\emph{S.~Sternberg} Lectures on Differential Geometry, Prentice Hall Publishers, 1964, 400p.
\bibitem{NT}
\emph{N.~Tanaka}, A differential-geometric study of strongly pseudoconvex manifolds, Lectures in Math., vol.9, Kyoto University, 1975
\bibitem{T}
\emph{S.~Tanno}, Variational problems on contact Riemannian manifolds, Transactions of the AMS, v.314, 1, 1989, p. 349--379
\bibitem{W}
\emph{S.~Webster} Pseudohermitian structures on a real hypersurface, J.Diff.Geom. 13(1978), p.25--31
\end{thebibliography}
\end{document}